\documentclass[11pt]{article}
\usepackage{amssymb, amsmath, amsthm, amsfonts,enumerate}
\usepackage{amsmath,setspace,scalefnt}
\usepackage{graphicx,tikz,caption,subcaption}
\usetikzlibrary{patterns}
\usetikzlibrary{shapes}
\usepackage{fullpage}
\usepackage{hyperref}
\usepackage{enumitem,verbatim}%,itemize
\usetikzlibrary{decorations.pathmorphing}
\usepackage{subcaption}

\tikzset{snake it/.style={decorate, decoration=snake}}
\usetikzlibrary{calc}

\definecolor{gray}{rgb}{0.25, 0.25, 0.25}

\newtheorem{theorem}{Theorem}[section]
\newtheorem{lemma}[theorem]{Lemma}
\newtheorem{cor}[theorem]{Corollary}

\theoremstyle{definition}

\theoremstyle{plain}

\theoremstyle{definition}

\theoremstyle{definition}

\theoremstyle{definition}

\theoremstyle{definition}
\newtheorem{defn}[theorem]{Definition}
\theoremstyle{definition}

\theoremstyle{definition}

\newcommand{\tk}{\mathsf{TK}}

\newcounter{propcounter}

\newcommand{\makenote}[2]
{

\smallskip

\noindent

\fbox{
\begin{minipage}{0.95\textwidth}

\def\temp{#1}
\ifx\temp\empty
\def\forlabel{$\bullet$}
\else
\def\forlabel{{\bfseries #1:}}
\fi

\begin{itemize}[label = \forlabel]
  #2/
\end{itemize}
\end{minipage}
}

\smallskip
}

\title{Rainbow clique subdivisions}

 \author{
Yan Wang
\footnote{School of Mathematical Sciences, CMA-Shanghai, Shanghai Jiao Tong University, Shanghai 200240, China (yan.w@sjtu.edu.cn)}
\thanks{Partially supported by National Key R\&D Program of China under Grant No. 2022YFA1006400, National Natural Science Foundation of China under grant No.12201400 and Explore X project of Shanghai Jiao Tong University}
}

\begin{document}
\maketitle

%\tableofcontents
%%%%%%%%%%%%%%%%%%%%%%%%%%%%%%%%%%%%%%%%%%%%%%%%
%%%%%%%%%%%%%%%%%%%%%%%%%%%%%%%%%%%%%%%%%%%%%%%%
%%%%%%%%%%%%%%%%%%%%%%%%%%%%%%%%%%%%%%%%%%%%%%%%
%%%%%%%%%%%%%%%%%%%%%%%%%%%%%%%%%%%%%%%%%%%%%%%%
%%%%%%%%%%%%%%%%%%%%%%%%%%%%%%%%%%%%%%%%%%%%%%%%
%\begin{itemize}
%	\item
%\end{itemize}

\begin{abstract} 
We show that for any integer $t \ge 2$, every properly edge colored $n$-vertex graph with average degree at least $(\log n)^{2+o(1)}$ contains a rainbow subdivision of a complete graph of size $t$.
Note this bound is within $(\log n)^{1+o(1)}$ factor of the lower bound. 
This also implies a result on the rainbow Tur\'{a}n number of cycles. 
\end{abstract}

\section{Introduction}\label{sec:intro}

Let $G$ be a graph. 
A \textit{subdivision} of $G$, denoted by $\mathsf{T}G$, is a graph obtained from $G$ by replacing each of its edges with internally vertex disjoint paths. 
%We call the vertices of $\mathsf{T}G$ corresponding to the vertices of $G$ its \textit{core} vertices.  
Subdivisions play an important role in graph theory.
One of the important results on subdivisions dates back to 1930s where Kuratowski \cite{Kur30} showed that a graph is not planar if and only if it contains a subdivision of a complete graph on five vertices or a subdivision of a complete bipartite graph with three vertices in each part. 

Mader \cite{Mad67} initiated the study of the relation between the average degree of a graph and the size of its largest clique subdivisions. 
For integer $t > 0$, let $d(t)$ be the minimum number $d$ such that every graph with average degree at least $d$ contains a subdivision of a complete graph $K_t$.
Mader \cite{Mad67} showed the existence of $d(t)$ in 1967.
Mader \cite{Mad67}, and independently Erd\H{o}s and Hajnal \cite{EH69} conjectured that $d(t) = O(t^2)$.
Subsequently, Mader \cite{Mad72} showed that $O(2^t)$ is an upper bound of $d(t)$.
In 1990s, Koml\'os and Szemer\'edi \cite{K-Sz-1, K-Sz-2}, and independently, Bollob\'as and Thomason \cite{BT98} confirmed this conjecture. 
As for lower bound, Jung \cite{Jung70} observed that disjoint union of complete regular bipartite graphs give the lower bound of $d(t) = \Omega(t^2)$.
Hence, $d(t) = \Theta(t^2)$.

In order to achieve a subdivision of a complete graph of size linear to the average degree, some additional conditions, such as minimum girth conditions, are needed to eliminate the extremal examples.
In fact, Mader \cite{Mad99} conjectured that every $C_4$-free graph of average degree $d$ contains a $\tk_{\Omega(d)}$.
K\"uhn and Osthus \cite{KO02, KO06} proved that every graph with sufficiently large girth contains a $\tk_{\delta(G)+1}$.
%subdivision of a complete graph of size larger than its minimum degree. 
They \cite{KO04} also showed the existence of $\tk_{d/\log^{12} d}$ in every $C_4$-free graph of average degree $d$.
In \cite{BLS15}, Balogh, Liu and Sharifzadeh proved Mader's conjecture assuming the graph is $C_6$-free.
Liu and Montgomery \cite{liu2017proof} completely resolved this conjecture recently. 
%Note that those proofs utilize the technique developed by  Koml\'os and Szemer\'edi \cite{K-Sz-1, K-Sz-2}.

For $\ell \in \mathbb{N}$, a \textit{balanced subdivision} of $G$, denoted by $\mathsf{T}G^{(\ell)}$, is a graph obtained from $G$ by replacing each of its edges with internally vertex disjoint paths of length exactly $\ell$. 
Thomassen \cite{Tho84, Tho85, thomassen} conjectured that for every constant $k \in \mathbb{N}$, there exists $d$ such that every graph with average degree at least $d$ contains a $\tk^{(\ell)}_k$ for some $\ell \in \mathbb{N}$.
Liu and Montgomery \cite{liu2020proof} confirmed Thomassen's conjecture. 
More recently, the author \cite{Wan23} showed that in every graph with average degree at least $d$  there is a $\tk_{\Omega(d^c)}^{(\ell)}$
for every constant $0 < c < 1/2$, which is improved to $\tk_{\Omega(\sqrt{d})}^{(\ell)}$ by Gil Fernández, Hyde, Liu, Pikhurko and Wu \cite{FHLPW22}, and Luan, Tang, Wang and Yang \cite{LTWY22} independently. 
Note that $\ell$ is a polylogarithmic function of the number of vertices of the graph in these results. 
Balanced clique subdivisions have also been studied extensively when restricting $\ell$ to be constant (see \cite{AKS03, ES74, Erd71,FS11,Jan21,Jia11,JS12,KP88,Tom22}).

In a graph with a proper edge coloring, we say that a subgraph is rainbow if all the edges have distinct colors. 
A rainbow variant of clique subdivision problems was considered by Jiang, Methuku and Yepremyan \cite{JMY21}.
They proved that every properly edge-colored graph on $n$ vertices with average degree at least $e^{c\sqrt{\log n}}$ contains a rainbow $\tk_t$.
Later, this upper bound on average degree was improved to $(\log n)^{60}$ by Jiang, Letzter, Methuku and Yepremyan \cite{JLMY21}.
Recently, Tomon \cite{Tom22} showed that $(\log n)^{6+o(1)}$ suffices. 
In this paper, we prove the following.

\begin{theorem}
\label{main_clique}
Let $t > 0$ be an integer.
Suppose $G$ is a properly edge colored graph on $n$ vertices with average degree at least $(\log n)^{2+o(1)}$.
Then $G$ contains a rainbow $\tk_t$.
\end{theorem}

%\begin{theorem}
%\label{main}
%Let $t > 0$ be an integer, $\varepsilon > 0$ and $n > 0$ be a sufficiently large integer. 
%Let $G$ be a proper edge colored graph on $n$ vertices.
%Suppose $d(G) \ge (\log n)^{1+\varepsilon}$.
%Then $G$ contains a rainbow $\tk_t$.
%\end{theorem}

We would like to point out that Theorem \ref{main_clique} is closely related to the study of rainbow Tur\'{a}n number.
Let $H$ be a graph. 
The Tur\'{a}n number $ex(n,H)$ is the maximum number of edges that a graph on $n$ vertices without a copy of $H$ can have.
Keevash, Mubayi, Sudakov and Verstra\"{e}te \cite{KMSV07} first introduced the following rainbow variant of Tur\'{a}n number. 
The rainbow Tur\'{a}n number $ex^*(n,H)$ is the maximum number of edges that a properly edge colored graph on $n$ vertices without a rainbow copy of $H$ can have.
In \cite{KMSV07}, Keevash, Mubayi, Sudakov and Verstra\"{e}te showed that $ex^*(n,H) = (1+o(1))ex(n,H)$ for non-bipartite $H$, thus determined the asymptotic value of $ex^*(n,H)$ by Erd\H{o}s-Stone-Simonovits Theorem \cite{ES66, ES46}. 
When $H$ is bipartite, determining $ex^*(n,H)$ is harder. 
In particular, much attention has been drawn on the study of $ex^*(n,C_{2k})$ where $C_{2k}$ is a cycle of length $2k$ (see \cite{DLS13, Jan20, KMSV07}), and Janzer \cite{Jan20} determined $ex^*(n,C_{2k}) = \Theta(n^{1+1/k})$.

It is well known that a graph with $n$ vertices without a cycle contains at most $n-1$ edges. 
It is natural to ask how many edges a properly edge colored graph on $n$ vertices without a rainbow cycle can have.
Equivalently, let $\mathcal{C}$ be the set of all cycles, it is interesting to determine $ex^*(n,\mathcal{C})$.
Keevash, Mubayi, Sudakov and Verstra\"{e}te \cite{KMSV07} showed that $ex^*(n,\mathcal{C}) = O(n^{4/3})$ and ask what number it should be.
Later, Das, Lee, Sudakov \cite{DLS13} improved the bound to $ne^{(\log n)^{1/2+o(1)}}$.
This was further improved by Janzer \cite{Jan20} to $O(n(\log n)^4)$.
The best upper bound is $n(\log n)^{2+o(1)}$ obtained recently by Tomon \cite{Tom22}.
It is easy to see that Theorem \ref{main_clique} implies $ex^*(n,\mathcal{C}) \le n (\log n)^{2+o(1)}$ (for example, take $t=3$).

\begin{cor}
\label{main_cycle}
Suppose $G$ is a properly edge colored graph on $n$ vertices with average degree at least $(\log n)^{2+o(1)}$.
Then $G$ contains a rainbow cycle.
\end{cor}

We remark that both Theorem \ref{main_clique} and Corollary \ref{main_cycle} are within $(\log n)^{1+o(1)}$ factor of the lower bound because of the following example due to Keevash, Mubayi, Sudakov and Verstra\"{e}te \cite{KMSV07}.
Consider $d$-dimensional hypercube $Q_d$: the vertices of $Q_d$ are all the subsets of $\{1,2,\ldots,d\}$ and the edges of $Q_d$ consist of all pairs of subsets of $[d]$ whose Hamming distance is exactly $1$.
Let $f$ be a proper edge coloring of $Q_d$ such that $f( \{X,X\setminus \{i\} \}) = i$ for $X \subseteq [d]$ and $i \in X$.
One can check $Q_d$ with such edge coloring $f$ contains no rainbow cycle. 
Moreover, the average degree of $Q_d$ is $d = \log n$.
This implies $ex^*(n,\mathcal{C}) \ge nd/2 = \Omega( n \log n)$.

\medskip

%The organization of this paper is as follows. 
%In Section \ref{sec:pre}, we introduce notations and lemmas that will be used in the proof. 
%We give a proof of Theorem \ref{main_clique} in Section \ref{sec:clique}.

The proof of Theorem \ref{main_clique} adopts the idea in \cite{Tom22} together with some new ideas.
First we generalize the definition of $\alpha$-maximal graphs to $\omega$-maximal graphs (see Definition \ref{def:omega}). 
We show that log-maximal graphs have good expansion property even after sampling the colors (see Lemma \ref{lem:expand}).
Using the sprinkling technique introduced in \cite{Tom22}, we show that from every vertex in a log-maximal graph one can reach more than half of the vertices via a rainbow path of logarithmic length avoiding a given set of vertices and colors (see Lemma \ref{lem:reachable}). 
Then it implies that any two vertices in a log-maximal graph can be connected by a rainbow path of small length upon removal of a set of vertices and colors of moderate size (see Lemma \ref{lem:diamter}). 
Finally we complete the proof by a greedy argument.

\subsection{Notations} \label{sec:notation}

For an integer $n \ge 1$, let $[n] = \{1,2,\ldots,n\}$.
Let $G$ be a graph. 
Let $V(G)$ and $E(G)$ be vertex set and edge set of $G$ respectively.
We define $v(G) = |V(G)|$ and $e(G) = |E(G)|$.
Let $X \subseteq V(G)$, we write $G - X$ for the induced subgraph of $G[V(G) \backslash X]$.
Let $X, Y \subseteq V(G)$, we write $G[X,Y]$ for the induced bipartite subgraph of $G$ with parts $X$ and $Y$.
We define $e_G(X,Y) = |E(G[X,Y])|$.
Let $d(G), \delta(G), \Delta(G)$ be the average degree, minimum degree and maximum degree of $G$ respectively.
For $v \in V(G)$, let $d_G(v)$ denote the degree of $v$ in $G$.
For $X \subseteq V(G)$, denote $N_G(X)$ the (external) neighborhood of $X$ in $G - X$. 
We omit the subscript if there is no confusion.
We also omit the floors and ceilings when they are not crucial.
All logarithms are base $2$. 

%For two vertices $u,v \in V(G)$, a $u,v$-path is a path with end vertices $u$ and $v$. 
%We use $\ell(P)$ to denote the length of a path $P$.
%The distance between two set of vertices $U,V$ in a graph $G$ is the minimum length of a $u,v$-path in $G$ with $u \in U$ and $v \in V$.

%For integer $i \ge 0$, we define the $i$-th ball around $X$ in $G$ to be the set of vertices that have distance at most $i$ from $X$ in $G$, denoted by  $B^i_G(X)$.
%For convenience, $B_G(X) = B^1_G(X) = X \cup N_G(X)$.

%In this section we show that for any integer $t \ge 2$, every proper edge colored $n$-vertex graph with average degree at least $O((\log n)^{1+o(1)})$ contains a rainbow $\tk_t$. 

\section{Preliminaries}
\label{sec:pre}

We need the following definitions. 

\begin{defn}
Let $G$ be a graph with proper edge coloring $f: E(G) \rightarrow R$.
Let $\phi: V(G) \rightarrow 2^{V(G) \cup R} $ be a mapping that assigns a set of (forbidden) vertices and colors for each vertex in $G$. 
For $X \subseteq V(G)$ and $Q \subseteq R$,
the \textit{restricted external neighborhood of $X$ in $G$ with respect to the colors $Q$} is 
$$N_{Q,\phi}(X) := \{ y \in V(G) \backslash X : \exists x \in X, xy \in E(G), f(xy) \in Q \backslash \phi(x), y \not\in \phi(x) \}$$
\end{defn}

\begin{defn}
Let $G$ be a graph with proper edge coloring $f: E(G) \rightarrow R$.
A \textit{rainbow $Q$-path} in $G$ is a path $v_1 v_2 \cdots v_k$ in $G$ such that $f(v_i v_{i+1}) \in Q$ for $i \in [k-1]$ and all $f(v_i v_{i+1})$ are distinct. 
\end{defn}

We need Markov's inequality (see \cite{probmethod}). 

\begin{lemma}\emph{(Markov's inequality)}
\label{markov}
Let $Y$ be a nonnegative random variable and $\alpha > 0$. Then we have 
$$\mathbb{P}(Y \ge \alpha \mathbb{E}[Y] ) \le \frac{1}{\alpha}.$$
\end{lemma}

We need the multiplicative Chernoff bound (see \cite{probmethod}). 

\begin{lemma}\emph{(Multiplicative Chernoff bound)}
\label{chernoff}
Let $X_1, \ldots, X_n$ be independent random variables taking values from $\{0,1\}$.
Let $X = \sum_{i=1}^n X_i$ and $\mu = \mathbb{E}[X]$.
Then we have
$$\mathbb{P}(X \le \frac{\mu}{2} ) \le e^{-\mu^2/8}.$$
\end{lemma}

We also need the following lemma, which is Lemma 2.4 in \cite{Tom22}.

\begin{lemma}
\label{concentration-tom22}
Let $p, p_c \in (0, 1]$, and $\lambda > 1$. Let $G$ be a bipartite graph with vertex classes $A$ and $B$, and
let $f : E(G) \rightarrow R$ be a proper edge coloring. Let $U \subseteq A$ be a random sample of vertices, each vertex
included independently with probability $p$, and let $Q \subseteq R$ be a random sample of colors, each included
independently with probability $p_c$. Let $\mu := E(|N_Q(U)|)$, and suppose that every vertex in $A$ has degree
at most $K$. If $K + |A| \le \frac{\mu}{32 \lambda \log(\lambda (p p_c)^{-1})}$, then
$$\mathbb{P} \left( |N_Q(U)| \le \frac{\mu}{64 \lambda \log(\lambda (p p_c)^{-1})} \right) \le 2e^{-\lambda}.$$
\end{lemma}

\section{Rainbow clique subdivisions}
\label{sec:clique}

In this section, we prove Theorem \ref{main_clique}.

\subsection{$\omega$-maximal graphs}

We generalize the definition of $\alpha$-maximal graphs in \cite{Tom22} as follows.

\begin{defn}
\label{def:omega}
Let $\omega: \mathbb{R}^+ \rightarrow \mathbb{R}^+$ be a function.
A graph $G$ is called \textit{$\omega$-maximal} if for every subgraph $H$ of $G$, we have
$$\frac{d(H)}{\omega(v(H))} \le \frac{d(G)}{\omega(v(G))}.$$
\end{defn}

It is easy to see that if $\omega(x) = x^{\alpha}$ then $\omega$-maximal graphs and $\alpha$-maximal graphs defined in \cite{Tom22} are the same.
Using the definition, it is not hard to see that an $\omega$-maximal graph $G$ has minimum degree at least $d(G)/2$. (It would also follow from Lemma 2.2 in \cite{JMY21})

\begin{lemma}
\label{lem:maximal}
Let $\omega: \mathbb{R}^+ \rightarrow \mathbb{R}^+$ be an increasing function.
Let $G$ be an $\omega$-maximal graph. 
Then $\delta(G) \ge d(G)/2$.
\end{lemma}

\begin{comment}
\begin{proof}
Let $v \in V(G)$ be a vertex of degree $\delta(G)$ and $H = G - \{v\}$. %$H = G[V(G) \backslash \{v\} ]$.
Since $G$ is $\omega$-maximal, we have $\frac{d(H)}{\omega(v(H))} \le \frac{d(G)}{\omega(v(G))}$, which is 
 $$\frac{d(G)v(G)-2\delta(G)}{(v(G) - 1) \omega(v(G) - 1) } \le \frac{d(G)}{\omega(v(G))}.$$
So we have
$$\delta(G) \ge \frac{1}{2} d(G) \left( v(G) - \frac{(v(G)-1) \omega(v(G)-1)}{\omega (v(G))} \right) \ge \frac{1}{2}d(G).$$
\end{proof}
\end{comment}

In the rest of this section, we take $\omega = x \mapsto \log x$.

\subsection{Expansions in log-maximal graphs}

We show that log-maximal graphs have good expansion property even after sampling the colors. 
The following lemma is similar to Lemma 2.6 in \cite{Tom22}.
However, we do not sample the vertices, which makes the proof simpler. 
Moreover, log-maximality is a bit more efficient than $\alpha$-maximality.
We include all the technical details as the calculations are affected by the choice of log-maximality.

\begin{lemma}
\label{lem:expand}
Let $0 < p_c \le 1$, $\lambda > 10^{8}$ and $n > 0$ be a sufficiently large integer such that $p_c \ge 1/\log n$. 
Let $G$ be a graph on $n$ vertices with proper edge coloring $f: E(G) \rightarrow R$ and $B \subseteq V(G)$ satisfying the following: 
\begin{itemize}
\item[(i)] $G$ is log-maximal;
\item[(ii)] $d(G) \ge \lambda^2 p_c^{-1} \log n \cdot \log (\lambda^{1/2} p_c^{-1} )$;
\item[(iii)] $\phi: V(G) \rightarrow 2^{V(G) \cup R}$ such that $|\phi(v)| \le \frac{d(G)}{8 \log n} \log(\frac{2n}{3|B|})$ for all $v \in V(G)$; 
\item[(iv)] $ 2 \le |B| \le \frac{n}{2}$.
\end{itemize}
Let $Q \subseteq R$ be a random subset of colors such that each color is chosen with probability $p_c$ independently.
Then with probability at least $1 - e^{-\Omega(\lambda^{1/2})}$, we have 
$$|N_{Q,\phi}(B)| \ge \min \left( \frac{|B|}{4}, 
 \frac{|B|\log(\frac{2n}{3|B|})}{8\log|B|}
\right).$$
\end{lemma}

\begin{proof}
Let $A = N_G(B)$.
Let $H$ be a bipartite graph with vertex classes $A$ and $B$ and edge set
$E(H) = \{xy: x \in B, y \in A \backslash \phi(x), f(xy) \in R \backslash \phi(x) \}$.
Let $H_Q$ be the subgraph of $H$ whose edges are colored with a color from $Q$, i.e.
$V(H_Q) = A \cup B$ and 
$E(H_Q) = \{xy: x \in B, y \in A \backslash \phi(x), f(xy) \in Q \backslash \phi(x) \}$.

Let $\Delta = \lambda^{1/2} p_c^{-1}$. 
Let $S = \{v \in A: |N_H(v)| \ge \Delta \}$ and $T = A \backslash S$.
We distinguish cases by the number of edges between $B$ and $T$ in $G$.

\medskip

\noindent \textbf{Case 1.} $e_G(B,T) \le \frac{d(G) |B|}{4 \log n} \log(\frac{2n}{3|B|})$.

First, we claim that $|S| \ge \min \{ \frac{|B|}{2}, 
 \frac{|B|}{4\log|B|} \log(\frac{2n}{3|B|})
\}.$
Otherwise, suppose $|S| < |B|/2$.
Let $C = V(G) \backslash B$.
Since $E(G) = E(G[B \cup S]) \cup E(G[C]) \cup E(G[B, T])$, %\cup E(G[B_1, B_0])
we have
\begin{equation*}
\begin{split}
{d(G[B \cup S])(|B| + |S|)}/{2} 
&=e(G[B \cup S]) \\
&\ge e(G) - e(G[C]) - e_G(B, T)  \\     %- e_G(B_1, B_0)
&= {d(G)n}/{2} - {d(G[C])|C|}/{2} - e_G(B, T)    %- |B_1||B_0|
\end{split}
\end{equation*}
As $G$ is log-maximal, we have
$$\frac{d(G[B \cup S])}{\log (|B| + |S|)} \le \frac{d(G)}{\log n},$$
and
$$\frac{d(G[C])}{\log |C|} \le \frac{d(G)}{\log n}.$$
Hence,
$$\frac{d(G)(|B| + |S|) \log(|B| + |S|) }{2\log n} \ge \frac{d(G)n}{2} - \frac{d(G)|C| \log |C|}{2\log n} - e_G(B, T).$$ % - |B_1||B_0|
Since  $e_G(B,T) \le \frac{d(G) |B|}{4 \log n} \log(\frac{2n}{3|B|})$,
\begin{equation*}
\begin{split}
(|B| + |S|) \log(|B| + |S|) 
&\ge n \log n - |C| \log |C| - \frac{2\log n}{d(G)} e_G(B, T) \\
&\ge |B| \log n + |C| (\log n - \log |C|) - \frac{2\log n}{d(G)} \frac{d(G) |B|}{4 \log n} \log(\frac{2n}{3|B|}) \\
&\ge |B| \log n - \frac{|B|}{2} \log(\frac{2n}{3|B|})
\end{split}
\end{equation*}
Since $|S| < |B|/2$,
\begin{equation*}
\begin{split}
|S| \log(\frac{3|B|}{2}) 
&\ge |B| (\log n - \log(\frac{3|B|}{2}) ) - \frac{|B|}{2} \log(\frac{2n}{3|B|}) = \frac{|B|}{2} \log(\frac{2n}{3|B|}).
\end{split}
\end{equation*}
Therefore,
$$|S|  
\ge \frac{|B|}{2\log(\frac{3|B|}{2})} \log(\frac{2n}{3|B|}) 
> \frac{|B|}{4\log|B|} \log(\frac{2n}{3|B|}).$$
This completes the proof the claim.

\medskip

Now let $W = N_{H_Q}(B) \cap S$.
For every vertex $y \in S$, we have 
$$\mathbb{P}(y \in W) = 1 - (1-p_c)^{|N_H(y) \cap B|} \ge 1 - (1-p_c)^{\Delta} = 1 - (1-p_c)^{\lambda^{1/2} p_c^{-1}} \ge 1 - e^{-\lambda^{1/2}}.$$
Thus, $\mathbb{E}[|W|] \ge |S|(1 - e^{-\lambda^{1/2}})$ and thus $\mathbb{E}[|S \setminus W|] \le |S| e^{-\lambda^{1/2}}$.
By Lemma \ref{markov} with $(Y,\alpha)_{\ref{markov}} = (|S \setminus W|, \frac{|S|}{2\mathbb{E}[|S\setminus W|]})$, we have
$$\mathbb{P}(| W| \le |S|/2)=\mathbb{P}(|S\setminus W| \ge |S|/2)=\mathbb{P}(|S\setminus W| \ge \frac{|S|}{2\mathbb{E}[|S\setminus W|]} \mathbb{E}[|S\setminus W|] ) \le \frac{2\mathbb{E}[|S\setminus W|]}{|S|} \le 2e^{-\lambda^{1/2}}.$$
Hence, we have
$$\mathbb{P}\left(|W| \ge \min \left( \frac{|B|}{4}, 
 \frac{|B|}{8\log|B|} \log(\frac{2n}{3|B|})
\right) \right) \ge \mathbb{P}(|W| \ge |S|/2) \ge 1 - 2 e^{-\lambda^{1/2}}.$$
This concludes the proof of Case 1.

\medskip

\noindent \textbf{Case 2.} $e_G(B,T) > \frac{d(G) |B|}{4 \log n} \log(\frac{2n}{3|B|})$.

Since $|\phi(v)| \le \frac{d(G)}{8 \log n} \log(\frac{2n}{3|B|})$ for all $v \in V(G)$ by (iii),
we have 
$$e_H(B,T) \ge e_G(B,T) - \sum_{v \in B} |\phi(v)| > \frac{d(G) |B|}{4 \log n} \log(\frac{2n}{3|B|}) - \frac{d(G) |B|}{8 \log n} \log(\frac{2n}{3|B|}) = \frac{d(G) |B|}{8 \log n} \log(\frac{2n}{3|B|}).$$

%So we can bound $|T|$ from below as 
%$$|T| \ge \frac{e_H(B_1,T)}{\Delta} \ge \frac{d |B_1|}{16\Delta \log n} \ge \frac{|B_1|}{4} \frac{d p_c }{4\lambda^{1/2}\log n}.$$
Let $T_1 = \{v \in T: |N_H(v)| < p_c^{-1}/2 \}$ and $T_2 = T \backslash T_1$.
If $e_H(B,T_1) \ge e_H(B,T_2)$, then 
let $W = N_{H_Q}(B) \cap T_1$. 
For every vertex $y \in T_1$, we have
$$\mathbb{P}(y \in W) = 1 - (1-p_c)^{|N_H(y) \cap B|} > 1 - (1-\frac{p_c|N_H(y) \cap B|}{2})=\frac{p_c|N_H(y) \cap B|}{2}.$$
Then 
$$\mathbb{E}|W| = \sum_{y \in T_1} \mathbb{P}(y \in W) > \frac{p_c}{2} \sum_{y \in T_1} |N_H(y) \cap B| = \frac{p_c}{2} e_H(B,T_1) \ge \frac{p_c}{4} e_H(B,T) >
\frac{d(G) p_c |B|}{32 \log n} \log(\frac{2n}{3|B|})
.$$

Otherwise, $e_H(B,T_1) < e_H(B,T_2)$.
Let $W = N_{H_Q}(B) \cap T_2$. 
For every vertex $y \in T_1$, we have
$$\mathbb{P}(y \in W) = 1 - (1-p_c)^{|N_H(y) \cap B|} > 1 - e^{-1/2} > 1/3.$$
Then 
$$\mathbb{E}|W| = \sum_{y \in T_2} \mathbb{P}(y \in W) > |T_2|/3 \ge  \frac{e_H(B,T_2)}{3\Delta} \ge  \frac{e_H(B,T)}{6\Delta} \ge \frac{d(G) p_c |B|}{48\lambda^{1/2} \log n} \log(\frac{2n}{3|B|}).$$

In both cases, let $\mathbb{E}[N_{H_Q}(B)] \ge \mathbb{E}|W| \ge \frac{d(G) p_c |B|}{48\lambda^{1/2} \log n} \log(\frac{2n}{3|B|}).$
We will apply Lemma \ref{concentration-tom22} with $(p,p_c,\lambda,G,A,B,K)_{\ref{concentration-tom22}}=(1,p_c,\lambda^{1/2},H[B,T],B,T,\Delta)$. Since $d(G) \ge \lambda^2 p_c^{-1} \log n \cdot \log (\lambda^{1/2} p_c^{-1} )$ and $p_c \ge 1/\log n$, we have
$$\frac{\mathbb{E}[N_{H_Q}(B)]}{32 \lambda^{1/2} \log(\lambda^{1/2} p_c^{-1})} \ge \frac{d(G) p_c |B| \log(\frac{2n}{3|B|})}{1536\lambda \log n \log(\lambda^{1/2} p_c^{-1})} \ge \frac{\lambda |B| \log(\frac{2n}{3|B|})}{1536} \ge \lambda^{1/2} p_c^{-1} + |B|.$$
Thus by Lemma \ref{concentration-tom22}, we have
$$\mathbb{P} \left( |N_{H_Q}(B)| \le \frac{\mathbb{E}[N_{H_Q}(B)]}{64 \lambda^{1/2} \log(\lambda^{1/2} p_c^{-1})} \right) \le 2e^{-\lambda^{1/2}}.$$

Therefore, with probability at least $1 - e^{-\Omega(\lambda^{1/2})}$, we have
$$|N_{H_Q}(B)| > \frac{\mathbb{E}[N_{H_Q}(B)]}{64 \lambda^{1/2} \log(\lambda^{1/2} p_c^{-1})} \ge \frac{\lambda |B| \log(\frac{2n}{3|B|})}{3072} \ge \min \left( \frac{|B|}{4}, 
 \frac{|B|\log(\frac{2n}{3|B|})}{8\log|B|}
\right)$$
since $d(G) \ge \lambda^2 p_c^{-1} \log n \cdot \log (\lambda^{1/2} p_c^{-1} )$ by (ii), $\lambda > 10^8$ and $p_c \ge 1/\log n$.
This completes the proof of Case 2.
\end{proof}

Now we show that every vertex in a log-maximal graph can reach many vertices by a rainbow path of moderate length. 
Note that the iterative neighbourhoods given by Lemma \ref{lem:expand} expand slightly faster due to the choice of log-maximality, which allows us to build a rainbow clique subdivision more efficiently than in \cite{Tom22}.

\begin{lemma}
\label{lem:reachable}
Let $0 < p_c \le 1$ and $n > 0$  be a sufficiently large integer such that $p_c \ge 1/\log n$. 
Let $G$ be a graph on $n$ vertices with a proper edge-coloring $f: E(G) \rightarrow R$ satisfying the following: 
\begin{itemize}
\item[(i)] $G$ is log-maximal;
\item[(ii)] $d := d(G) \ge \lambda^3 p_c^{-1} (\log n)^2$ where $\lambda \ge (\log \log n)^{10}$; % and $C$ is a large constant;
\item[(iii)]  $\phi_0 \subseteq V(G) \cup R$ with $|\phi_0| \le {d}/{16 \log n}$.
\end{itemize}
Let $Q \subseteq R$ be a random subset of colors such that each color is chosen with probability $p_c$ independently.
Then for every $v \in V(G)$, with probability more than $1/2$, more than $n/2$ vertices of $G$ can be reached from $v$ by a rainbow $(Q \backslash \phi_0)$-path of length $O(\log n \cdot \log \log n)$ avoiding (forbidden) vertices in $\phi_0$.
%with edges in colors $Q \backslash \phi_0$
\end{lemma}

\begin{proof}
Let $l = 32 \log n \cdot \log \log n $; so $l \le {d}/{32 \log n}$.
We adopt the ``sprinkling" technique. 
We sample colors with different probability in each round so that the final distribution of colors is the same after this process ends.
More precisely, we define $q_i$ for $i \in [l]$ as follows: 
$q_1 = p_c / 2$ and $q_i = q$ for $i \in [l] \backslash \{1\}$ where $1-p_c = (1-q_1)(1-q)^{l-1}$. 
Thus, $q = \Theta(p_c/l) = \Theta(p_c / (\log n \cdot \log \log n))$.

Fix $v \in V(G)$.
Let $Q_i$ be a random sample of $Q$ such that each color is chosen with probability $q_i$ independently for $i \in [l]$.
We define $\phi_i$, $S_i$ and $B_i$ recursively as follows.
$$\phi_0(v) := \phi_0, \phi_0(x) := \emptyset, \forall x \in V(G) \backslash \{v\}, $$
$$S_1 = B_1 := \{x \in N(v) \backslash \phi_0(v) : f(xv) \in Q_1 \backslash \phi_0(v) \}.$$
%% The following description is not accurate, because our definition $B_i$ requires the color of $i$-th segment is from $Q_i$. 
%For $i \in [l]$, $B_i$ is a subset of the vertices that are reachable from $v$ by a rainbow $((Q_1 \cup \cdots \cup Q_i) \backslash \phi_0)$-path of length at most $i$ avoiding vertices in $\phi_{0}$.
For each $x \in B_i$, let $P_{vx}^i$ be an arbitrary rainbow $((Q_1 \cup \cdots \cup Q_i) \backslash \phi_0)$-path from $v$ to $x$ of length at most $i$.
For $i \in [l-1]$, we define $\phi_i(x)$ to be the union of $\phi_0$ and vertices and colors used in $P_{vx}^i$ for $x \in B_i$; and $\phi_i(x) = \emptyset$ for $x \not\in B_i$.
For $i \in [l-1]$, define
$$S_{i+1} := \{ y \in N(B_i) \backslash (B_i \cup \phi_0): \exists x \in B_i, y \not\in \phi_i(x), xy \in E(G), f(xy) \in Q_{i+1} \backslash \phi_i(x) \}$$
$$B_{i+1} := B_i \cup S_{i+1}.$$

We want to apply Lemma \ref{lem:expand} to $(G,B,p_c,\lambda)_{\ref{lem:expand}} = (G,B_i,q_i,\lambda)$ for every $i \in [l]$, and show that $B_i$ expands substantially. 
Note that 
$$d \ge \lambda^3 p_c^{-1} (\log n)^2 > \lambda^2  q_i^{-1} \log n \cdot \log (\lambda^{1/2} q_i^{-1})$$ 
and 
$$|\phi_i(x)| \le 2l + |\phi_0| \le 2 \cdot \frac{d }{32 \log n} + \frac{d}{16 \log n} = \frac{d}{8 \log n}.$$
Moreover, as $G$ is log-maximal, $\delta(G) \ge d/2$ by Lemma \ref{lem:maximal}.
We have 
$$\mathbb{E}[|B_1|] \ge \delta(G) q_1 - |\phi_0| \ge \frac{dp_c}{4} - \frac{d}{16 \log n} > 8.$$ 
Hence, by Lemma \ref{chernoff} with $(X)_{\ref{markov}} = (|B_1|)$, we have
$$\mathbb{P}(|B_1| \le \frac{\mathbb{E}[|B_1|]}{2} ) \le e^{-\mathbb{E}[|B_1|]^2/8} < \frac{1}{4}.$$
So with probability at least $3/4$, $|B_1| \ge {\mathbb{E}[|B_1|]}/{2} > 2$, thus $|B_i| \ge |B_1| > 2$ for all $i \in [l]$.

%Since $|B_i|$ is non-decreasing, it suffices to show that $|B_1| \ge 2 \lambda q_i^{-1}$.
%Indeed, as $G$ is $\alpha$-maximal, $\delta(G) \ge d/2$.
%So $\mathbb{E}|B_1| \ge \delta(G) q_0 = d p_c/4.$
%By Lemma \ref{chernoff}, we have 
%$$\mathbb{P}(|B_1| < d p_c / 8) < e^{-\frac{(\mathbb{E}|B_1|)^2}{8}} < e^{-\frac{d^2 p_c^2}{128}} = e^{-\Omega(\lambda)}.$$
%So with probability $1-e^{-\Omega(\lambda)}$, $|B_1| \ge d p_c / 8 \ge \lambda^2 \alpha^{-1} p_c^{-1} \log n > 2 \lambda q_i^{-1}$.

Therefore, by Lemma \ref{lem:expand}  with probability at least $1 - e^{-\Omega(\lambda^{1/2})}$, we have 
$$|S_{i+1}| \ge |N_{Q_{i+1},\phi_i}(B_i)| 
\ge \min \left( \frac{|B_i|}{4}, 
 \frac{|B_i|\log(\frac{2n}{3|B_i|})}{8\log(|B_i|)}
\right).$$ 
With probability $3/4 - l \cdot e^{-\Omega(\lambda^{1/2})} > 1/2$, this is true for all $i \in [l-1]$.
If $|B_i| < ({2n}/{3})^{1/3}$, then $|S_{i+1}| \ge {|B_i|}/{4}$.
Otherwise, $|S_{i+1}| \ge  {|B_i|\log(\frac{2n}{3|B_i|})}/{8\log(|B_i|)}$.

Let $r_1 \ge 1$ be the minimum integer such that $|B_{r_1}| \ge (\frac{2n}{3})^{1/3} $.
It is easy to see that $r_1 \le {\log(2n/3)}/{3\log(5/4)} = O(\log n)$.

For $r_1 \le i \le l$, let $\delta_i > 0$ be such that $|B_i| = (2n/3)^{1-\delta_i}$.
Now let $r_1 \le r_2 \le l$ be the minimum integer such that $|B_{r_2}| > {n}/{2} =  ({2n}/{3})^{1 - \log_{2n/3}(4/3)} $.
Thus, for $r_1 \le i < r_2$,
$$|B_{i+1}| = |B_i| + |S_{i+1}| \ge |B_i| \left(1 + \frac{\delta_i \log(2n/3)}{8 (1-\delta_i) \log(2n/3)} \right) = |B_i| \left( 1 + \frac{\delta_i}{8(1-\delta_i)} \right) \ge |B_i| \frac{1}{1-\delta_i/8}.$$
Hence,
$$1 - \delta_{i+1} \ge 1 - \delta_i + \log_{2n/3} \left( \frac{1}{1-\delta_i/8} \right).$$
$$\delta_{i+1} \le \delta_i + \log_{2n/3} \left( 1-\delta_i/8 \right) = \delta_i + \frac{\log(1-\delta_i/8)}{\log (2n/3)} \le \delta_i (1 - \frac{1}{8 \log (2n/3)}).$$

By definition of $r_2$, we have $\delta_{r_2-1} \ge \log_{2n/3}(4/3)$.
Therefore, $r_2 \le r_1 + 10 \log n \cdot \log \log n \le l$.
So $|B_l| > {n}/{2}$, which concludes the proof.
\end{proof}

As a corollary, we are able to show that every two vertices in a log-maximal graph on $n$ vertices can be connected by a path of length at most $O(\log n \cdot \log \log n)$ upon forbidding a moderate size of vertices and colors.
In other word, its small diameter property is robust.

\begin{lemma}
\label{lem:diamter}
Let $n > 0$  be a sufficiently large integer. 
Let $G$ be a graph on $n$ vertices with proper edge coloring $f: E(G) \rightarrow R$ satifying the following: 
\begin{itemize}
\item[(i)] $G$ is log-maximal;
\item[(ii)] $d := d(G) \ge 4 \lambda^3 (\log n)^2$ where $\lambda \ge (\log \log n)^{10}$; % and $C$ is a large constant;
\item[(iii)]  $\phi_0 \subseteq V(G) \cup R$ with $|\phi_0| \le {d}/{16 \log n}$.
\end{itemize}
For any two vertices $u,v \in V(G)$, there exists a rainbow $(R \backslash \phi_0)$-path from $u$ to $v$ of length $O(\log n \cdot \log \log n)$ avoiding vertices in $\phi_0$.
\end{lemma}

\begin{proof}
Let $p_c = 1/2$ and $(R_u, R_v)$ be a partition of $R$ such that each color appears in $R_u$ with probability $1/2$.
One can view $R_u$ (resp. $R_v$) as a random subset of colors $R$ such that each color is chosen with probability $p_c$ independently.

Let $B_u$ (resp. $B_v$) be the vertices of $G$ can be reached from $u$ (resp. $v$) by a rainbow $(R_u \backslash \phi_0)$-path (resp. $(R_v \backslash \phi_0)$-path) of length $O(\log n \cdot \log \log n)$ avoiding (forbidden) vertices in $\phi_0$.
Apply Lemma \ref{lem:reachable} to $G,u$ with $p_c = 1/2$ (resp. to $G,v$ with $p_c = 1/2$), we have
with probability more than $1/2$, $|B_u| \ge n/2$ (resp. $|B_v| \ge n/2$).

Therefore, with positive probability, there exists a partition $(R_u, R_v)$ of $R$ such that $|B_u| > n/2$ and $|B_v| > n/2$, thus $B_u \cap B_v \ne \emptyset$.
Let $w \in B_u \cap B_v$ and $P_{uw}$ (resp. $P_{vw}$) be a rainbow $(R_u \backslash \phi_0)$-path (resp. $(R_v \backslash \phi_0)$-path) from $u$ (resp. $v$) to $w$ of length $O(\log n \cdot \log \log n)$ avoiding (forbidden) vertices in $\phi_0$.
Choose $w$ such that $|V(P_{uw})| + |V(P_{vw})|$ is minimum.
Therefore, $P_{uw} \cup P_{vw}$ is a desired rainbow path.
\end{proof}

\subsection{Proof of Theorem \ref{main_clique}}

\begin{proof}
We prove that for every $\varepsilon > 0$ and $n$ sufficiently large, if $G$ is a properly edge colored graph on $n$ vertices of average degree at least $(\log n)^{2+\varepsilon}$, then $G$ contains a rainbow copy of $\tk_t$.
By passing onto a subgraph, we may assume that $d(G) = (\log n)^{2+\varepsilon}$.
Suppose $f: E(G) \rightarrow R$ is the proper edge coloring of $G$.
Let $H$ be a log-maximal subgraph of $G$ and $m = v(H)$.
So 
$$d(H) \ge \frac{\log m}{\log n} d(G) = \log m \cdot (\log n)^{1+\varepsilon} \ge (\log m)^{2+\varepsilon}$$
 and $\delta(H) \ge d(H)/2 > (\log m)^{2+\varepsilon}/2$.
 Note that $m \ge d(H) \ge (\log n)^{1+\varepsilon}$; so $m$ is also sufficiently large.

Let $v_1,v_2,\cdots,v_t$ be $t$ distinct vertices in $H$.
% with pairwise distance at least $2$. Such vertices exist as we can pick vertices one by one greedily. 
Let $K \subseteq {t \choose 2}$ be a maximal collection of pairs such that there exists a family of pairwise internally disjoint rainbow paths $\mathcal{P} = \{k \in K: P_k\}$ such that
\stepcounter{propcounter}
\begin{enumerate}[label = {\bfseries \Alph{propcounter}\arabic{enumi}}]
  \item For each $\{i,j\} \in K$, $P_{\{i,j\}}$ is a rainbow path of length $O(\log m \cdot \log \log m)$ from $v_i$ to $v_j$; \label{main-1}
  \item No colors appear more than once in $\{f(e): e \in P, P \in \mathcal{P} \}$. \label{main-2}
\end{enumerate}
If $K = {t \choose 2}$, then the graph formed by all the paths in $\mathcal{P}$ is a desired rainbow $\tk_t$.
Hence, we may assume that there exist distinct $i, j \in [t]$ such that $\mathcal{P}$ contains no such path from $v_i$ to $v_j$.

Let $\phi_0$ be the union of vertices and colors in the paths in $\mathcal{P}$ except $v_i$ and $v_j$.
Note that 
$|\phi_0| \le {t \choose 2} \cdot O(\log m \cdot \log \log m) < {d(H)}/{16 \log m}$
and 
$d(H) \ge (\log m)^{2+\varepsilon} \ge 4 (\log \log m)^{30} (\log m)^2$.
Apply Lemma \ref{lem:diamter} with $(G,\lambda,\phi_0)_{\ref{lem:diamter}} = (H,(\log \log m)^{10},\phi_0)$, we obtain a rainbow $(R \backslash \phi_0)$-path $P_{\{i,j\}}$ from $v_i$ to $v_j$ of length $O( \log m \cdot \log \log m)$ avoiding vertices in $\phi_0$.
Hence, $K \cup \{\{i,j\}\}$ and $\mathcal{P} \cup \{ P_{\{i,j\}} \}$ contradict the maximality of $K$.
This completes the proof.
\end{proof}

\section*{Acknowledgements}

We thank the anonymous referees for their careful reading and suggestions.
We also thank Hong Liu for bringing our attention to this problem. 
We are grateful to Istvan Tomon for useful discussions.

%%%%%%%%%%%%%%%%%%%%%%%%%%%%%%%%%%%%%%%%%%%%%%%%%%%%%%%%%%%%%%%%%%%

\bibliographystyle{abbrv}
\bibliography{rainbow_clique_subdivision}

\end{document}